\documentclass[11pt]{amsart}

\usepackage{amsmath}
\usepackage{amssymb}
\usepackage{mathrsfs}

\newtheorem{theorem}{Theorem}[section]
\newtheorem{claim}[theorem]{Claim}

\newtheorem{lemma}[theorem]{Lemma}

\newtheorem{proposition}[theorem]{Proposition}
\newtheorem{corollary}[theorem]{Corollary}

\theoremstyle{definition}
\newtheorem{definition}[theorem]{Definition}

\newtheorem{question}[theorem]{Question}

\theoremstyle{remark}

\newcount\skewfactor
\def\mathunderaccent#1#2 {\let\theaccent#1\skewfactor#2
\mathpalette\putaccentunder}
\def\putaccentunder#1#2{\oalign{$#1#2$\crcr\hidewidth
\vbox to.2ex{\hbox{$#1\skew\skewfactor\theaccent{}$}\vss}\hidewidth}}
\def\name{\mathunderaccent\tilde-3 }

\def\smallbox#1{\leavevmode\thinspace\hbox{\vrule\vtop{\vbox
   {\hrule\kern1pt\hbox{\vphantom{\tt/}\thinspace{\tt#1}\thinspace}}
   \kern1pt\hrule}\vrule}\thinspace}

\newcommand{\cf}{{\rm cf}}

\def\qedref#1{$\qed_{\reforiginal{#1}}$}

\title{Polarized relations at singulars over successors}
\author{Shimon Garti}
\address{Einstein Institute of Mathematics,
 The Hebrew University of Jerusalem,
 Jerusalem 91904, Israel}
\email{shimon.garty@mail.huji.ac.il}


\begin{document}
\let\labeloriginal\label
\let\reforiginal\ref
\def\ref#1{\reforiginal{#1}}
\def\label#1{\labeloriginal{#1}}

\begin{abstract}
Erd\H{o}s, Hajnal and Rado asked whether $\binom{\aleph_{\omega_1}}{\aleph_2}\rightarrow \binom{\aleph_{\omega_1}}{\aleph_0}_2$ and whether $\binom{\aleph_{\omega_1}}{\aleph_2}\rightarrow \binom{\aleph_{\omega_1}}{\aleph_1}_2$.
We shall prove that both relations are independent over ZFC.
We shall also prove that $\binom{\mu}{\aleph_2}\rightarrow \binom{\mu}{\aleph_2}_2$ is independent over ZF for some $\mu>\cf(\mu)=\omega_1$.
\end{abstract}

\maketitle

\newpage

\section{Introduction}

Let $\kappa\leq\lambda$ be infinite cardinals.
We shall say that $\binom{\lambda}{\kappa}\rightarrow\binom{\alpha}{\beta}_\chi$ iff for every coloring $c:\lambda\times\kappa\rightarrow\chi$ one can find $A\subseteq\lambda,B\subseteq\kappa$ such that ${\rm otp}(A)=\alpha,{\rm otp}(B)=\beta$ and $c\upharpoonright(A\times B)$ is constant.
The case in which $\lambda=\kappa^+$ is of particular interest, and the strong form $\binom{\kappa^+}{\kappa}\rightarrow\binom{\kappa^+}{\kappa}_2$ has been investigated quite thoroughly.

It is still unknown whether the positive relation $\binom{\kappa^+}{\kappa}\rightarrow\binom{\kappa^+}{\kappa}_2$ is consistent when $\kappa$ is a successor cardinal.
The first case is $\kappa=\aleph_1$, and the possible consistency of $\binom{\aleph_2}{\aleph_1}\rightarrow\binom{\aleph_2}{\aleph_1}_2$ is an interesting open problem, a variant of which appeared in \cite[Question 27]{MR0280381}.
In this light, Erd\H{o}s, Hajnal and Rado asked whether a much weaker statement like $\binom{\aleph_2}{\aleph_1}\rightarrow\binom{\aleph_0}{\aleph_1}_2$ is consistent, see \cite[Problem 12]{MR0202613}.

Soon after the discovery of forcing, Prikry proved in \cite{MR0297577} that $\binom{\aleph_2}{\aleph_1}\nrightarrow\binom{\aleph_0}{\aleph_1}_2$ is consistent.
A few years later, Laver proved in \cite{MR673792} that $\binom{\aleph_2}{\aleph_1}\rightarrow\binom{\aleph_0}{\aleph_1}_2$ and even the much stronger relation $\binom{\aleph_2}{\aleph_1}\rightarrow\binom{\aleph_1}{\aleph_1}_\omega$ are consistent.
Hence \cite[Problem 12]{MR0202613} of Erd\H{o}s, Hajnal and Rado is independent of ZFC.

A related problem arises if one replaces $\aleph_1$ by a singular cardinal $\mu$ such that $\cf(\mu)=\omega_1$.
Thus one may wonder whether $\binom{\aleph_{\omega_1}}{\aleph_2}\rightarrow \binom{\aleph_{\omega_1}}{\aleph_0}_2$ and even whether $\binom{\aleph_{\omega_1}}{\aleph_2}\rightarrow \binom{\aleph_{\omega_1}}{\aleph_1}_2$ are provable, see \cite[Problem 15]{MR0202613}.
We shall prove that both relations are independent of ZFC.

The above problem is phrased with respect to $\aleph_{\omega_1}$, but the general question is about every $\mu>\cf(\mu)=\omega_1$.
The ultimate relation $\binom{\mu}{\aleph_2}\rightarrow \binom{\mu}{\aleph_2}_2$ is far-fetched, though we do not know whether it provably fails in ZFC.
However, we will show that this relation is independent over ZF by proving the positive direction under the assumption ${\rm AD}+V=L(\mathbb{R})$.
Determinacy will also help us to establish a positive answer to an old problem from \cite{MR0280381} about cube polarized relations.

Our notation is standard for the most part.
We shall follow \cite{MR795592} with respect to arrows notation.
We suggest \cite{MR2768693} for background in pcf theory and \cite{MR3075383} for background in polarized relations.
We employ the Jerusalem forcing notation, hence $p\leq q$ reads $p$ is weaker than $q$.

Finally, I would like to thank the referees of the paper for a careful reading of the original manuscript and a lot of helpful suggestions.

\newpage

\section{Between MA and GMA}

In this section we address the first part of \cite[Problem 15]{MR0202613} by proving the independence of $\binom{\aleph_{\omega_1}}{\aleph_2}\rightarrow \binom{\aleph_{\omega_1}}{\aleph_0}_2$.
Ahead of proving our statements let us recall that if $\kappa\geq\cf(\kappa)>\omega$ and one forces Martin's axiom with $2^\omega>\kappa$ then one obtains $\binom{\kappa}{\omega}\rightarrow\binom{\kappa}{\omega}_2$, see \cite{MR0371652}.
In particular, Martin's axiom with $2^\omega>\omega_1$ implies $\binom{\omega_1}{\omega}\rightarrow\binom{\omega_1}{\omega}_2$.
Therefore, a natural attempt to get the consistency of $\binom{\aleph_2}{\aleph_1}\rightarrow\binom{\aleph_2}{\aleph_1}_2$ (or at least something in this direction) would be one of the generalized versions of Martin's axiom to uncountable cardinals.

We shall see below that this attempt fails.
This will be done by proving $\binom{\aleph_2}{\aleph_1}\nrightarrow\binom{\aleph_0}{\aleph_1}_2$ under one of the traditional generalizations of Martin's axiom, and in the same context we will also have $\binom{\aleph_{\omega_1}}{\aleph_2}\nrightarrow \binom{\aleph_{\omega_1}}{\aleph_0}_2$.
The known generalizations are similar, and we shall use Shelah's version from \cite{MR0505492} dubbed henceforth as the generalized Martin's axiom.

\begin{theorem}
\label{thmshelahgma} Generalized Martin's Axiom. \newline
One can force $2^{\aleph_0}=\aleph_1\wedge 2^{\aleph_1}>\aleph_2$, and if $\mathbb{P}$ is a forcing notion of size less than $2^{\aleph_1}$ satisfying the following three requirements:
\begin{enumerate}
\item [$(a)$] Each pair of compatible conditions has a least upper bound in $\mathbb{P}$.
\item [$(b)$] Every countable increasing sequence of conditions has a least upper bound in $\mathbb{P}$.
\item [$(c)$] If $\{p_i:i<\aleph_2\}\subseteq\mathbb{P}$ then there is a club $C\subseteq\aleph_2$ and a regressive function $f:\aleph_2\rightarrow\aleph_2$ so that for $\alpha,\beta\in C\cap S^{\aleph_2}_{\aleph_1}$ if $f(\alpha)=f(\beta)$ then $p_\alpha\parallel p_\beta$.
\end{enumerate}
then there is a generic filter $G\subseteq\mathbb{P}$ which intersects any given collection of $\kappa$ dense subsets, when $\kappa<2^{\aleph_1}$.
If $\kappa$ satisfies $\alpha<\kappa\Rightarrow\alpha^{\aleph_0}<\kappa$ then the assumption $|\mathbb{P}|<2^{\aleph_1}$ can be omitted.
\end{theorem}

\hfill \qedref{thmshelahgma}

Our small component in the theorems below will be $\aleph_1$, but in most of the statements we can replace $\aleph_1$ by larger successor cardinals.
For the consistency of the negative direction $\binom{\aleph_{\omega_1}}{\aleph_2}\nrightarrow \binom{\aleph_{\omega_1}}{\aleph_0}_2$ we invoke a result of Prikry.
An early paper of Prikry from \cite{MR0297577} contains a forcing construction which proves the consistency of $\binom{\omega_2}{\omega_1}\nrightarrow\binom{\omega_0}{\omega_1}_2$.
The theorem generalizes to $\binom{\kappa}{\omega_1}\nrightarrow\binom{\omega_0}{\omega_1}_2$ for every regular uncountable $\kappa$.
The proof is exactly as in \cite{MR0297577}, and for completeness we unfold it upon replacing $\omega_2$ by $\kappa$.
Let us begin with the following:

\begin{definition}
\label{defpmatrix} Prikry matrix. \newline
Suppose that $\kappa=\cf(\kappa)\geq\aleph_2$. \newline
A $(\kappa\times\omega_1)$-Prikry matrix is a family $(A_{\alpha,\eta}:\alpha\in\kappa,\eta\in\omega_1)$ of subsets of $\omega_1$ satisfying the following properties:
\begin{enumerate}
\item [$(\aleph)$] For every $\alpha\in\kappa$ and $\zeta<\eta<\omega_1, A_{\alpha,\zeta}\cap A_{\alpha,\eta}=\varnothing$.
\item [$(\beth)$] For every $\alpha\in\kappa, \bigcup\{A_{\alpha,\eta}:\eta\in\omega_1\}=\omega_1$.
\item [$(\gimel)$] For every $\{\alpha_n:n\in\omega\}\subseteq\kappa$ and every sequence $(\eta_n:n\in\omega)\subseteq\omega_1, |\omega_1-\bigcup\{A_{\alpha_n,\eta_n}:n\in\omega\}|\leq\aleph_0$.
\end{enumerate}
\end{definition}

We use the set notation $\{\alpha_n:n\in\omega\}$ to express the fact that $m\neq n\Rightarrow\alpha_m\neq\alpha_n$, and the sequence notation $(\eta_n:n\in\omega)$ to allow repetitions.

\begin{claim}
\label{clmneg} If there is a $(\kappa\times\omega_1)$-Prikry matrix then $\binom{\kappa}{\aleph_1}\nrightarrow\binom{\aleph_0}{\aleph_1}_2$.
\end{claim}

\par\noindent\emph{Proof}. \newline
Suppose that $(A_{\alpha,\eta}:\alpha\in\kappa,\eta\in\omega_1)$ is a $(\kappa\times\omega_1)$-Prikry matrix.
Define $c:\kappa\times\omega_1\rightarrow 2$ by $c(\alpha,\beta)=0$ iff $\beta\in A_{\alpha,\eta}$ and $\eta$ is an even ordinal.
Assume that $A\in[\kappa]^{\aleph_0}$ and $B\in[\omega_1]^{\aleph_1}$.
Enumerate the elements of $A$ by $\{\alpha_n:n\in\omega\}$.
Fix an even ordinal $\eta\in\omega_1$ and an odd ordinal $\zeta\in\omega_1$.
Since $|\omega_1-\bigcup\{A_{\alpha_n,\eta_n}:n\in\omega\}|\leq\aleph_0$ and $|B|=\aleph_1$ we can choose an ordinal $\beta\in B\cap\bigcup_{n\in\omega}A_{\alpha_n,\eta}$, so $\beta\in A_{\alpha_n,\eta}$ for some $n\in\omega$.
It follows that $c(\alpha_n,\beta)=0$.
Similarly one can choose an ordinal $\gamma\in B\cap\bigcup_{n\in\omega}A_{\alpha_n,\zeta}$, so $\gamma\in A_{\alpha_m,\eta}$ for some $m\in\omega$.
It follows that $c(\alpha_m,\gamma)=1$.
We conclude, therefore, that $A\times B$ is $c$-polychromatic as desired.

\hfill \qedref{clmneg}

Our next goal is to force a $(\kappa\times\omega_1)$-Prikry matrix.
The forcing is identical with that of \cite{MR0297577}, upon replacing $\aleph_2$ by $\kappa$.

\begin{definition}
\label{defpmforcing} Prikry-matrix forcing. \newline
Let $\kappa=\cf(\kappa)\geq\aleph_2$. \newline
A condition $p\in\mathbb{P}$ is a pair $(S,F)=(S^p,F^p)$ such that:
\begin{enumerate}
\item [$(i)$] $S\subseteq\kappa\times\omega_1\times\omega_1, |S|\leq\aleph_0$.
\item [$(ii)$] For every $\alpha\in\kappa,\beta\in\omega_1$ there is at most one ordinal $\eta\in\omega_1$ such that $(\alpha,\eta,\beta)\in S$.
\item [$(iii)$] $F$ is a function from ${\rm dom}(F)$ into $\omega_1$, where ${\rm dom}(F)$ is countable and the elements of ${\rm dom}(F)$ are functions $f:{\rm dom}(f)\rightarrow\omega_1$ so that ${\rm dom}(f)\subseteq\kappa$ and $|{\rm dom}(f)|\leq\aleph_0$.
\item [$(iv)$] For every $\beta\in\omega_1$ if there are $\alpha\in\kappa,\eta\in\omega_1$ such that $(\alpha,\eta,\beta)\in S$ then for every $f\in{\rm dom}(F)$ with $F(f)\leq\beta$ there exists $\alpha'\in{\rm dom}(f)$ such that $(\alpha',f(\alpha'),\beta)\in S$.
\end{enumerate}
If $p,q\in\mathbb{P}$ then $p\leq q$ iff $S^p\subseteq S^q$ and $F^p\subseteq F^q$.
\end{definition}

The following lemma shows that forcing with $\mathbb{P}$ preserves cardinals.
Moreover, it implies that the forcing notion $\mathbb{P}$ satisfies stronger properties which are parallel to the requirements of Theorem \ref{thmshelahgma}.

\begin{lemma}
\label{lemknaster} Let $\kappa=\cf(\kappa)\geq\aleph_2$ and let $\mathbb{P}$ be the associated Prikry-matrix forcing.
\begin{enumerate}
\item [$(\aleph)$] If $p,q\in\mathbb{P}$ and $p\parallel q$ then $r=(S^p\cup S^q,F^p\cup F^q)$ is a least upper bound for $p$ and $q$.
\item [$(\beth)$] Likewise, if $(p_n:n\in\omega)$ is $\leq_{\mathbb{P}}$-increasing then its union is a least upper bound. In particular, $\mathbb{P}$ is $\aleph_1$-complete.
\item [$(\gimel)$] $\mathbb{P}$ is $\kappa$-Knaster.
\end{enumerate}
\end{lemma}

\par\noindent\emph{Proof}. \newline
The first two parts follow directly from the definition, and the last part follows from a Delta-system argument.

\hfill \qedref{lemknaster}

Observe that for $\aleph_2$ we conclude from the above lemma that the Prikry-matrix forcing satisfies the requirements in Theorem \ref{thmshelahgma}.
We shall define now some dense subsets of $\mathbb{P}$.
For every pair $(\alpha,\beta)\in\kappa\times\omega_1$ let $D_{\alpha,\beta}=\{p\in\mathbb{P}: \exists\eta\in\omega_1,(\alpha,\eta,\beta)\in S^p\}$.
For every function $f$ which satisfies $(iii)$ of Definition \ref{defpmforcing} we let $E_f=\{p\in\mathbb{P}:f\in{\rm dom}(F^p)\}$.

\begin{lemma}
\label{lemdensity} Let $\kappa=\cf(\kappa)\geq\aleph_2$ and let $\mathbb{P}$ be the associated Prikry-matrix forcing. \newline
Then $D_{\alpha,\beta}$ is a dense open subset of $\mathbb{P}$ for every $(\alpha,\beta)\in\kappa\times\omega_1$, and $E_f$ is a dense open subset of $\mathbb{P}$ for every $f$ which satisfies $(iii)$ of Definition \ref{defpmforcing}.
\end{lemma}

\par\noindent\emph{Proof}. \newline
Directly from the definitions.

\hfill \qedref{lemdensity}

We can prove now the following generalization of \cite{MR0297577}:

\begin{theorem}
\label{thmmt2} Assume that $\kappa=\cf(\kappa)\geq\aleph_2$ and $\theta^\omega=\theta$ for every $\theta<\kappa$ of uncountable cofinality. \newline
Then one can force the existence a $(\kappa\times\omega_1)$-Prikry matrix.
\end{theorem}

\par\noindent\emph{Proof}. \newline
Let $\mathbb{P}$ be the Prikry-matrix forcing for $\kappa$.
Since $\mathbb{P}$ is $\aleph_1$-complete, it adds no $\omega$-subsets and hence $\theta^\omega=\theta$ for every $\theta<\kappa$ of uncountable cofinality in the generic extension.
Let $G\subseteq\mathbb{P}$ be $V$-generic.
We shall argue that there exists a $(\kappa\times\omega_1)$-Prikry matrix in $V[G]$.

Let us describe the sets in our Prikry matrix.
For every $\alpha\in\kappa$ and every $\eta\in\omega_1$ define:
$$
A_{\alpha,\eta}=\{\beta\in\omega_1: \exists S \exists F, (S,F)\in G \wedge (\alpha,\eta,\beta)\in S\}.
$$
We claim that $\{A_{\alpha,\eta}:\alpha\in\kappa,\eta\in\omega_1\}$ is a Prikry matrix.

Clearly, $A_{\alpha,\eta}\subseteq\omega_1$ for every $\alpha\in\kappa,\eta\in\omega_1$.
Suppose that $\alpha\in\kappa$ and $\zeta<\eta<\omega_1$.
Assume towards contradiction that $A_{\alpha,\zeta}\cap A_{\alpha,\eta}\neq\varnothing$ and choose an ordinal $\beta\in A_{\alpha,\zeta}\cap A_{\alpha,\eta}$.
Choose $p,q\in G$ such that $p\Vdash\check{\beta}\in\name{A}_{\alpha,\zeta}$ and $q\Vdash\check{\beta}\in\name{A}_{\alpha,\eta}$.
Since $G$ is directed we can find $r\in G$ such that $p,q\leq r$.
We conclude that $r\Vdash\check{\beta}\in\name{A}_{\alpha\zeta}\cap \name{A}_{\alpha,\eta}$ and hence $(\alpha,\zeta,\beta),(\alpha,\eta,\beta)\in S^r$.
This is impossible, however, since $\zeta\neq\eta$ and by virtue of Definition \ref{defpmforcing}$(ii)$.

Fix $\alpha\in\kappa$ and $\beta\in\omega_1$.
By Lemma \ref{lemdensity} we can choose a condition $p\in G\cap D_{\alpha,\beta}$.
This means that $(\alpha,\eta,\beta)\in S^p$ and hence $p\Vdash\check{\beta}\in\name{A}_{\alpha,\eta}$ for some $\eta\in\omega_1$.
We conclude, therefore, that $V[G]\models\bigcup_{\eta\in\omega_1}A_{\alpha,\eta}=\omega_1$.

Finally, let $\name{A}$ be a name of an element of $[\kappa]^{\aleph_0}$ and let $\name{\eta}$ be a name of an $\omega$-sequence of ordinals from $\omega_1$.
Fix a condition $p$ which forces these facts.
We may assume that $p$ forces that $A=\{\alpha_n:n\in\omega\}$ and $\eta=\langle\eta_n:n\in\omega\rangle$, since $\mathbb{P}$ is $\aleph_1$-complete so one can form an $\omega$-increasing sequence of conditions above $p$, each of which forces a value to another element of $A$ and $\eta$, and then take an upper bound.

Our goal is to find a condition $r\in G$ and an ordinal $\gamma\in\omega_1$ such that $r\Vdash\forall\beta\geq\gamma, \beta\in\bigcup\{A_{\alpha_n,\eta_n}:n\in\omega\}$.
To do this, we define a function $f$ as follows.
We let ${\rm dom}(f)=\{\alpha_n:n\in\omega\}$ and $f(\alpha_n)=\eta_n$ for every $n\in\omega$.
We choose a condition $q\in G\cap E_f$, so $f\in{\rm dom}(F^q)$, and we let $\gamma=F^q(f)$.

Suppose that $\gamma\leq\beta\in\omega_1$.
Choose an ordinal $\alpha\in\kappa$ and extend $q$ to a condition $r\in G\cap D_{\alpha,\beta}$, so $(\alpha,\eta,\beta)\in S^q$.
By Definition \ref{defpmforcing}$(iv)$ we see that $(\alpha_n,f(\alpha_n),\beta)\in S^q$ for some $n\in\omega$.
We see that $r\Vdash \check{\beta}\in\name{A}_{\alpha_n,f(\alpha_n)}=\name{A}_{\alpha_n,\eta_n}$, so we are done.

\hfill \qedref{thmmt2}

Based on Prikry's result, we can show that the first part of \cite[Problem 15]{MR0202613} is independent of ZFC:

\begin{theorem}
\label{thmehr15} The polarized relation $\binom{\aleph_{\omega_1}}{\aleph_2}\rightarrow \binom{\aleph_{\omega_1}}{\aleph_0}$ is independent of ZFC.
\end{theorem}

\par\noindent\emph{Proof}. \newline
The positive direction of $\binom{\aleph_{\omega_1}}{\aleph_2}\rightarrow \binom{\aleph_{\omega_1}}{\aleph_0}$ follows from the stronger relation $\binom{\aleph_{\omega_1}}{\aleph_0}\rightarrow \binom{\aleph_{\omega_1}}{\aleph_0}$ by monotonicity.
This stronger relation holds if one forces Martin's axiom with $2^\omega>\aleph_{\omega_1}$.

For the opposite direction, let us show that $\binom{\aleph_{\omega_1}}{\aleph_2}\nrightarrow \binom{\aleph_{\omega_1}}{\aleph_0}$ in Prikry's model.
So we force as in \cite{MR0297577} and we choose a coloring $c:\omega_2\times\omega_1\rightarrow 2$ which exemplifies the negative relation $\binom{\aleph_2}{\aleph_1}\nrightarrow\binom{\aleph_0}{\aleph_1}_2$.
Fix an increasing and continuous sequence of cardinals $\langle\mu_\gamma:\gamma\in\omega_1\rangle$ such that $\omega_{\omega_1}=\bigcup_{\gamma\in\omega_1}\mu_\gamma$.
For every ordinal $\alpha\in\omega_{\omega_1}$ let $\gamma(\alpha)$ be the unique ordinal in $\omega_1$ such that $\mu_{\gamma(\alpha)}\leq\alpha<\mu_{\gamma(\alpha)+1}$.
Given $\alpha\in\omega_{\omega_1}$ and $\beta\in\omega_2$ we define $d(\alpha,\beta)=c(\beta,\gamma(\alpha))$.

Assume that $A\subseteq\omega_{\omega_1},|A|=\aleph_{\omega_1}$ and $B\subseteq\omega_2,|B|=\aleph_0$.
We claim that $d\upharpoonright(A\times B)$ is not constant.
To see this notice that the set $I=\{\gamma\in\omega_1:\exists\alpha\in A,\gamma(\alpha)=\gamma\}$ is of size $\aleph_1$ since $A$ is unbounded in $\omega_{\omega_1}$.
By the choice of the coloring $c$ we know that $c\upharpoonright(B\times I)$ is not constant.

Pick up two pairs $(\beta_0,\gamma_0),(\beta_1,\gamma_1)\in B\times I$ such that $c(\beta_0,\gamma_0)=0$ and $c(\beta_1,\gamma_1)=1$.
Choose $\alpha_0,\alpha_1\in A$ such that $\gamma(\alpha_0)=\gamma_0$ and $\gamma(\alpha_1)=\gamma_1$.
By definition, $d(\alpha_0,\beta_0)=c(\beta_0,\gamma_0)=0$ and $d(\alpha_1,\beta_1)=c(\beta_1,\gamma_1)=1$, so the proof is accomplished.

\hfill \qedref{thmehr15}

As noted by one of the referees, a more general statement stands behind the above argument, and we phrase it for completeness:

\begin{proposition}
\label{propreferee} Let $\kappa,\lambda$ and $\mu$ be cardinals. \newline
Then $\binom{\kappa}{\lambda}\rightarrow\binom{\kappa}{\mu}_2$ implies $\binom{\cf(\kappa)}{\lambda}\rightarrow\binom{\cf(\kappa)}{\mu}_2$.
\end{proposition}

\par\noindent\emph{Proof}. \newline
Assume that $c:\cf(\kappa)\times\lambda\rightarrow 2$ is a coloring.
Choose a continuous increasing sequence of cardinals $(\kappa_i:i\in\cf(\kappa))$ so that $\kappa=\bigcup_{i\in\cf(\kappa)}\kappa_i$.
For every $\alpha\in\kappa$ there is a unique ordinal $i\in\cf(\kappa)$ such that $\kappa_i\leq\alpha<\kappa_{i+1}$, call it $i(\alpha)$.
Define a coloring $d:\kappa\times\lambda\rightarrow 2$ by letting $d(\alpha,\beta)=c(i(\alpha),\beta)$.
Since $\binom{\kappa}{\lambda}\rightarrow\binom{\kappa}{\mu}_2$ one can find $A\in[\kappa]^\kappa,B\in[\lambda]^\mu$ so that $d\upharpoonright(A\times B)$ is constantly $\ell$ for some $\ell\in\{0,1\}$.

Let $H=\{i(\alpha):\alpha\in A\}$ and notice that $|H|=\cf(\kappa)$ as $A$ is unbounded in $\kappa$.
If $\gamma\in H,\beta\in B$ then $\gamma=i(\alpha)$ for some $\alpha\in A$ and then $c(\gamma,\beta)=c(i(\alpha),\beta)=d(\alpha,\beta)=\ell$.
Therefore, the product $H\times B$ exemplifies the relation $\binom{\cf(\kappa)}{\lambda}\rightarrow\binom{\cf(\kappa)}{\mu}_2$ with respect to the coloring $c$.
Since $c$ was arbitrary, we are done.

\hfill \qedref{propreferee}

It is interesting to compare statements over $\omega$ which can be forced by Martin's axiom with parallel statements over $\omega_1$ under the generalized Martin's axiom or Baumgartner's axiom.
A systematic study in this direction is carried out in \cite{MR1278025}, and many applications of Martin's axiom can be forced over higher cardinals.
The above theorem shows, however, that not everything is possible:

\begin{corollary}
\label{corgma} If one forces $2^\omega=\omega_1$ and the generalized Martin's axiom with $2^{\omega_1}>\omega_2$ then one obtains the negative relation $\binom{\aleph_2}{\aleph_1}\nrightarrow\binom{\aleph_0}{\aleph_1}_2$ and hence $\binom{\aleph_2}{\aleph_1}\nrightarrow\binom{\aleph_2}{\aleph_1}_2$.
\end{corollary}

\par\noindent\emph{Proof}. \newline
Let $\mathbb{P}$ be the Prikry-matrix forcing at $\omega_2$.
As mentioned before, it follows from Lemma \ref{lemknaster} that $\mathbb{P}$ satisfies the requirements of Theorem \ref{thmshelahgma}.
Hence by forcing the generalized Martin's axiom with $2^{\omega_1}>\omega_2$ we add a generic $G\subseteq\mathbb{P}$ which intersects every prescribed collection of $\omega_2$ dense sets.
This is sufficient for an $(\omega_2\times\omega_1)$-Prikry matrix, so we are done.

\hfill \qedref{corgma}

We make the comment that strong positive relations of the form $\binom{\aleph_2}{\aleph_1}\rightarrow_{\mathcal{A}} \binom{\aleph_2}{\aleph_1}_2$ hold under the generalized Martin's axiom when we take a certain collection of colorings determined by $\mathcal{A}$, see \cite[Theorem 2.5]{MR3610266}.

It is not clear whether a positive relation with finite sets at the large component is consistent with the generalized Martin's axiom.
The positive relation $\binom{\aleph_2}{\aleph_1}\rightarrow \binom{2}{\aleph_1}_\omega$ has a considerable consistency strength, as proved by Donder and Levinski in \cite{MR1024901}.
Hence $\binom{\aleph_2}{\aleph_1}\nrightarrow \binom{2}{\aleph_1}_\omega$ is consistent with the generalized Martin's axiom since it can be forced without large cardinals.
Donder and Levinski showed that $\binom{\aleph_2}{\aleph_1}\rightarrow \binom{2}{\aleph_1}_\omega$ implies $\binom{\aleph_2}{\aleph_1}\rightarrow \binom{n}{\aleph_1}_\omega$ for every $n\in\omega$, and a proof can be found in \cite{MR1034564}.
This invites for the following question:

\begin{question}
\label{qdonderlevinski} Is it consistent that $\binom{\aleph_2}{\aleph_1}\rightarrow \binom{2}{\aleph_1}_\omega$ but $\binom{\aleph_2}{\aleph_1}\nrightarrow \binom{\aleph_0}{\aleph_1}_\omega$?
\end{question}

\newpage

\section{Scales and saturated ideals}

In this section we address the second part of \cite[Problem 15]{MR0202613} by showing that the relation $\binom{\aleph_{\omega_1}}{\aleph_2}\rightarrow \binom{\aleph_{\omega_1}}{\aleph_1}_2$ is independent.
The negative relation $\binom{\aleph_{\omega_1}}{\aleph_2}\nrightarrow \binom{\aleph_{\omega_1}}{\aleph_1}_2$ holds, for example, if one forces the generalized Martin's axiom with $2^{\omega_1}>\omega_2$ using the results of the previous section and monotonicity.
For the consistency of the positive direction we shall need the following:

\begin{definition}
\label{defscattered} Scattered families. \newline
Assume that $\mu>\cf(\mu)=\theta$ and let $\kappa=\mu^+$. \newline
A family of sets $\mathcal{A}=\{A_\alpha:\alpha\in\kappa\}$ will be called $\mu$-scattered iff the following requirements are met:
\begin{enumerate}
\item [$(\alpha)$] $A_\alpha\subseteq\mu$ for every $\alpha\in\kappa$.
\item [$(\beta)$] $|A_\alpha|=\theta$ for every $\alpha\in\kappa$.
\item [$(\gamma)$] If $I\in[\kappa]^\kappa$ and $\{a_\alpha:\alpha\in I\}$ satisfies $a_\alpha\subseteq A_\alpha,|a_\alpha|=\theta$ for every $\alpha\in I$ then $|\bigcup\{a_\alpha:\alpha\in I\}|=\mu$.
\end{enumerate}
\end{definition}

A scattered family will prove helpful when we try to lift polarized relations over a regular cardinal to the parallel relation with a singular cardinal sharing the same cofinality.

\begin{claim}
\label{clmscattered} Assume that:
\begin{enumerate}
\item [$(a)$] $\mu>\cf(\mu)=\theta, \kappa=\mu^+, 2^\theta<\kappa$.
\item [$(b)$] $\mathcal{A}=\{A_\alpha:\alpha\in\kappa\}$ is $\mu$-scattered.
\item [$(c)$] $\chi<\theta$ and $\binom{\theta^+}{\theta}\rightarrow\binom{\theta}{\theta}_\chi$.
\end{enumerate}
Then $\binom{\mu}{\theta^+}\rightarrow\binom{\mu}{\theta}_\chi$. \newline
If one assumes the stronger relation $\binom{\theta^+}{\theta}\rightarrow\binom{\theta^+}{\theta}_\chi$ and $2^{\theta^+}<\kappa$ then one obtains the stronger conclusion $\binom{\mu}{\theta^+}\rightarrow\binom{\mu}{\theta^+}_\chi$.
\end{claim}

\par\noindent\emph{Proof}. \newline
Let $c:\mu\times\theta^+\rightarrow\chi$.
For every $\alpha\in\kappa$ let $c_\alpha=c\upharpoonright(A_\alpha\times\theta^+)$.
Since $\binom{\theta^+}{\theta}\rightarrow\binom{\theta}{\theta}_\chi$ which can also be written as $\binom{\theta}{\theta^+}\rightarrow\binom{\theta}{\theta}_\chi$ we see that for every $\alpha\in\kappa$ there are $a_\alpha\in[A_\alpha]^\theta$ and $b_\alpha\in[\theta^+]^\theta$ such that $c_\alpha\upharpoonright(a_\alpha\times b_\alpha)$ is homogeneous with the color $i_\alpha\in\chi$.

Since $\kappa=\cf(\kappa)>2^\theta$ and $\chi<\theta$ we see that there are $I\in[\kappa]^\kappa, b\in[\theta^+]^\theta$ and a fixed $i\in\chi$ such that $\alpha\in I\Rightarrow b_\alpha=b\wedge i_\alpha=i$.
Let $a=\bigcup\{a_\alpha:\alpha\in I\}$ and notice that $|a|=\mu$ since $\mathcal{A}$ is $\mu$-scattered.

We claim that $c\upharpoonright(a\times b)$ is constantly $i$.
To see this, suppose that $\sigma\in a$ and $\tau\in b$.
Choose $\alpha\in\kappa$ such that $\sigma\in a_\alpha$.
It follows that $c(\sigma,\tau)=c_\alpha(\sigma,\tau)=i_\alpha=i$ and hence $\binom{\mu}{\theta^+}\rightarrow\binom{\mu}{\theta}_\chi$.
Moreover, if one assumes that the stronger relation $\binom{\theta^+}{\theta}\rightarrow\binom{\theta^+}{\theta}_\chi$ holds then the same argument yields $\binom{\mu}{\theta^+}\rightarrow\binom{\mu}{\theta^+}_\chi$, but we need $2^{\theta^+}<\kappa$ in order to unify the $b_\alpha$s into the same $b$, thus the proof is accomplished.

\hfill \qedref{clmscattered}

Our next objective is to show that scattered families exist.
We need a few basic concepts from pcf theory.
Suppose that $\mu>\cf(\mu)=\theta$ and $(\mu_\gamma:\gamma\in\theta)$ is an increasing sequence of regular cardinals such that $\mu=\bigcup_{\gamma\in\theta}\mu_\gamma$.
If $f,g\in\prod_{\gamma\in\theta}\mu_\gamma$ and $J$ is an ideal over $\theta$ then we shall say that $f<_Jg$ iff $\{\gamma\in\theta:g(\gamma)\leq f(\gamma)\}\in J$.
Usually, our ideal $J$ will be $J^{\rm bd}_\theta$, the ideal of bounded subsets of $\theta$.

A scale $(f_\alpha:\alpha\in\kappa)$ in the product $(\prod_{\gamma\in\theta}\mu_\gamma,J)$ is a sequence of elements of $\prod_{\gamma\in\theta}\mu_\gamma$ such that $\alpha<\beta\Rightarrow f_\alpha<_Jf_\beta$ and for every $h\in\prod_{\gamma\in\theta}\mu_\gamma$ there is an ordinal $\alpha\in\kappa$ such that $h<_Jf_\alpha$.
By a fundamental theorem of Shelah's pcf theory, if $\mu>\cf(\mu)=\theta$ and $\kappa=\mu^+$ then there exists an increasing sequence $(\mu_\gamma:\gamma\in\theta)$ of regular cardinals such that $\mu=\bigcup_{\gamma\in\theta}\mu_\gamma$ and there is a scale $(f_\alpha:\alpha\in\kappa)$ in $(\prod_{\gamma\in\theta}\mu_\gamma,J^{\rm bd}_\theta)$.
A proof can be found in \cite{MR1318912} or in \cite[Theorems 2.23 and 2.26]{MR2768693}.

\begin{proposition}
\label{propscattexist} If $\mu$ is a singular cardinal then there exists a $\mu$-scattered family.
\end{proposition}

\par\noindent\emph{Proof}. \newline
Let $\theta=\cf(\mu)$ and $\kappa=\mu^+$.
Choose an increasing sequence $(\mu_\gamma:\gamma\in\theta)$ of regular cardinals such that $\mu=\bigcup_{\gamma\in\theta}\mu_\gamma$.
Let $(f_\alpha:\alpha\in\kappa)$ be a scale in $(\prod_{\gamma\in\theta}\mu_\gamma,J^{\rm bd}_\theta)$.
For every $\alpha\in\kappa$ let $A_\alpha={\rm rang}(f_\alpha)$, so $A_\alpha\subseteq\mu$ for every $\alpha\in\kappa$.
Without loss of generality the range of each $f_\alpha$ is unbounded in $\mu$ and hence $|A_\alpha|=\theta$ for every $\alpha\in\kappa$.

Fix a set $I\in[\kappa]^\kappa$ and a collection $\{a_\alpha:\alpha\in I\}$ so that $a_\alpha\in[A_\alpha]^\theta$ for every $\alpha\in I$.
Let $a=\bigcup\{a_\alpha:\alpha\in I\}$.
If there is a set $y\in[\theta]^\theta$ such that $|a\cap\mu_\gamma|=\mu_\gamma$ for every $\gamma\in y$ then $|a|=|\bigcup_{\gamma\in y}\mu_\gamma|=\mu$ and we are done.
By way of contradiction assume that no such $y$ exists.
Hence there is an ordinal $\gamma_0\in\theta$ such that $\gamma\in(\gamma_0,\theta)\Rightarrow \eta_\gamma=\sup(a\cap\mu_\gamma)<\mu_\gamma$.
Define $h\in\prod_{\gamma\in\theta}\mu_\gamma$ by $h(\gamma)=0$ when $\gamma\leq\gamma_0$ and $h(\gamma)=\eta_\gamma$ when $\gamma>\gamma_0$.

Choose an ordinal $\alpha\in I$ such that $h<_{J^{\rm bd}_\theta}f_\alpha$.
For some $\delta_0\in\theta$ we see that if $\delta_0\leq\delta<\theta$ then $h(\delta)<f_\alpha(\delta)$.
Since $A_\alpha={\rm rang}(f_\alpha)$ and $a_\alpha$ is unbounded in $A_\alpha$ we can choose a sufficiently large $\delta$ so that $\gamma_0<\delta<\theta$ and $h(\delta)<f_\alpha(\delta)\in a_\alpha$.
However, $\sup(a_\alpha\cap\mu_\delta)\leq\sup(a\cap\mu_\delta) = \eta_\delta=h(\delta)<f_\alpha(\delta)$, a contradiction.

\hfill \qedref{propscattexist}

Recall that if Martin's axiom holds and $2^\omega>\omega_1$ then $\binom{\omega_1}{\omega}\rightarrow\binom{\omega_1}{\omega}_2$.
If $\theta$ is supercompact then one can force $\mathfrak{s}_\theta>\theta^+$ and obtain $\binom{\theta^+}{\theta}\rightarrow\binom{\theta^+}{\theta}_2$.
For this result see \cite[Theorem 2.4]{MR3509813}, which is based on \cite[Claim 1.2]{MR3201820}.
If $\theta$ is a successor of a regular cardinal and there exists a huge cardinal above $\theta$ then one can force $\binom{\theta^+}{\theta}\rightarrow\binom{\theta}{\theta}_2$, see \cite[Theorem 3.2]{MR3610266}.
These facts justify the following:

\begin{proposition}
\label{corehr15} Assume that $\mu>\cf(\mu)=\theta$.
\begin{enumerate}
\item [$(\aleph)$] If $\theta=\aleph_0$ or $\theta$ is supercompact then one can force $\binom{\mu}{\theta^+}\rightarrow\binom{\mu}{\theta}_2$.
\item [$(\beth)$] If $\theta$ is a successor of a regular cardinal and there is a huge cardinal above $\theta$ then one can force $\binom{\mu}{\theta^+}\rightarrow\binom{\mu}{\theta}_2$.
\item [$(\gimel)$] In particular, if there is a huge cardinal then one can force the relation $\binom{\aleph_{\omega_1}}{\aleph_2}\rightarrow \binom{\aleph_{\omega_1}}{\aleph_1}_2$.
\end{enumerate}
\end{proposition}

\par\noindent\emph{Proof}. \newline
Begin with a singular cardinal $\mu>\cf(\mu)=\theta$ so that $2^\theta<\mu$ (or even $2^{\theta^+}<\mu$ for a stronger conclusion).
Such a cardinal satisfies part $(a)$ of Claim \ref{clmscattered}, and part $(b)$ follows from Proposition \ref{propscattexist}.
Part $(c)$ can be forced as described in the paragraph before the statement of this proposition.
Applying Claim \ref{clmscattered} we are done.

\hfill \qedref{corehr15}

Along with the results of the previous section we obtain a full answer to the original question of \cite{MR0202613}.
One may wonder, however, if the strongest relation $\binom{\aleph_{\omega_1}}{\aleph_2}\rightarrow \binom{\aleph_{\omega_1}}{\aleph_2}_2$ is provably false in ZFC.
We shall discuss such relations in the next section, but we drop the axiom of choice in order to obtain strong positive relations.

\newpage

\section{Strong relations without choice}

The possible consistency of $\binom{\aleph_2}{\aleph_1}\rightarrow \binom{\aleph_2}{\aleph_1}_2$ with ZFC is still open.
However this positive relation is provable under the axiom of determinacy, and more is provable under the additional assumption that $V=L(\mathbb{R})$.
Other variations of the polarized relation under AD appear in the literature, see e.g. \cite{MR2984052}.
We shall also use determinacy in order to address the cube problem, which formally is still open in ZFC.

\begin{theorem}
\label{thmadcube} Assume AD.
\begin{enumerate}
\item [$(\aleph)$] $\binom{\aleph_2}{\aleph_1}\rightarrow \binom{\aleph_2}{\aleph_1}_{2^{\aleph_0}}$.
\item [$(\beth)$] $\left( \begin{smallmatrix} \aleph_2 \\ \aleph_1 \\ \aleph_0 \end{smallmatrix} \right) \rightarrow \left( \begin{smallmatrix} \aleph_2 \\ \aleph_1 \\ \aleph_0 \end{smallmatrix} \right)_n$ for every $n\in\omega$.
\item [$(\gimel)$] If one assumes ${\rm AD}+V=L(\mathbb{R})$ then
$\binom{\mu}{\aleph_2}\rightarrow \binom{\mu}{\aleph_2}_{2^{\aleph_0}}$ for some $\mu>\cf(\mu)=\omega_1$.
\end{enumerate}
\end{theorem}

\par\noindent\emph{Proof}. \newline
Recall that $\aleph_1$ and $\aleph_2$ are measurable under AD, as proved first by Solovay.
Likewise, ${\rm AC}_{\aleph_0}$ holds under AD, and we shall use these facts below.
We shall prove a general claim which will be helpful for all three statements.
Fix an $\aleph_1$-complete ultrafilter $\mathscr{U}$ over $\aleph_1$ and an $\aleph_2$-complete ultrafilter $\mathscr{V}$ over $\aleph_2$.
The claim is that if we are given a sequence $(c_n:n\in\omega)$ where $c_n:\omega_2\times\omega_1\rightarrow 2$ for every $n\in\omega$ then there are $A\in\mathscr{V}, B\in\mathscr{U}$ so that $c_n\upharpoonright(A\times B)$ is constant for every $n\in\omega$ simultaneously.

Suffice it to prove that for every $n\in\omega$ there are $A_n\in\mathscr{V}, B_n\in\mathscr{U}$ so that $c_n\upharpoonright(A_n\times B_n)$ is constantly $i_n$ for some $i_n\in\{0,1\}$.
Indeed, if we prove this statement then we can define $A=\bigcap\{A_n:n\in\omega\}$ and $B=\bigcap\{B_n:n\in\omega\}$.
Thus we have $A\in\mathscr{V}, B\in\mathscr{U}$ by the completeness of these ultrafilters.
But now we see that $c_n\upharpoonright(A\times B)$ is constantly $i_n$ for every $n\in\omega$ simultaneously.
Remark that we use here ${\rm AC}_{\aleph_0}$, as we must choose $A_n$ and $B_n$, but ${\rm AC}_{\aleph_0}$ is at our disposal under AD.

So fix $n\in\omega$ and focus on the coloring $c_n$.
For every $\beta\in\omega_1$ let $j_\beta\in\{0,1\}$ be such that $S_\beta^{j_\beta} = \{\alpha\in\omega_2:c_n(\alpha,\beta)=j_\beta\}\in\mathscr{V}$.
For some $B_n\in\mathscr{U}$ and a fixed $i_n\in\{0,1\}$ we will have $\beta\in B_n\Rightarrow j_\beta=i_n$.
Define $A_n=\bigcap\{S_\beta^{i_n}:\beta\in B_n\}$, and conclude that $A_n\in\mathscr{V}$ by $\aleph_2$-completeness.
Now $c_n\upharpoonright(A_n\times B_n)$ is constantly $i_n$ as desired.

We proceed to part $(\aleph)$ of our theorem.
We may identify $2^{\aleph_0}$ with the collection of $\omega$-sequences of $0$ and $1$, so assume that $c:\omega_2\times\omega_1\rightarrow{}^\omega 2$ is given.
For every $n\in\omega$ let $c_n:\omega_2\times\omega_1\rightarrow 2$ be the $n$th place of $c$, to wit $c_n(\alpha,\beta)$ is the $n$th place in the sequence $c(\alpha,\beta)$ for every $\alpha\in\omega_2,\beta\in\omega_1$.
Let $A\in\mathscr{V},B\in\mathscr{U}$ be such that $c_n''(A\times B)=\{i_n\}$ for every $n\in\omega$.
Here we use the general claim proved above.
Let $\eta=(i_n:n\in\omega)$, so $\eta\in{}^\omega 2$.
By definition, $c\upharpoonright(A\times B)$ is constantly $\eta$, so the first statement of the theorem is proved.

Part $(\beth)$ is basically the same.
For notational simplicity we focus on the case in which $n=2$, the general proof for every $n\in\omega$ is just the same.
We are given now a coloring $d:\omega_2\times\omega_1\times\omega\rightarrow 2$.
For every $n\in\omega$ we define $d_n:\omega_2\times\omega_1\rightarrow 2$ by letting $d_n(\alpha,\beta)=d(\alpha,\beta,n)$.
For every $n\in\omega$ there are $A_n\in\mathscr{V},B_n\in\mathscr{U}$ and $i_n\in\{0,1\}$ so that $c_n''(A_n\times B_n)=\{i_n\}$.
Let $C\in[\omega]^\omega$ be such that $n\in C\Rightarrow i_n=i$ for some fixed $i\in\{0,1\}$.
Finally, let $A=\bigcap_{n\in C}A_n, B=\bigcap_{n\in C}B_n$ so $A\in\mathscr{V}$ and $B\in\mathscr{U}$.
One can verify that $d\upharpoonright(A\times B\times C)$ is constantly $i$, so the second part has been established.

Finally, we wish to prove that $\binom{\mu}{\aleph_2}\rightarrow \binom{\mu}{\aleph_2}_{2^{\aleph_0}}$ under ${\rm AD}+V=L(\mathbb{R})$ for some $\mu>\cf(\mu)=\omega_1$.
Suffice it to prove that $\binom{\mu}{\aleph_2}\rightarrow \binom{\mu}{\aleph_2}_2$ since then we can use ${\rm AC}_\omega$ to get $2^{\aleph_0}$ colors by the general claim from the beginning of the proof.
We emphasize that the strengthening to $2^{\aleph_0}$ colors using the general claim is based on the fact that $\aleph_1$ is measurable, which is correct under ${\rm AD}+V=L(\mathbb{R})$.

Recall that under ${\rm AD}+V=L(\mathbb{R})$ one can prove that $\Theta$ is regular, and it is a limit of measurable cardinals.
In fact, if $\kappa<\Theta$ is regular then $\kappa$ is measurable.
Moreover, for every such $\kappa$ the filter generated by the $\omega$-closed unbounded subsets of $\kappa$ is a $\kappa$-complete ultrafilter over $\kappa$.
A proof of all these facts appears in \cite[Theorem 8.27]{MR2768698}.

Let $(\kappa_\alpha:\alpha\in\omega_1)$ be an increasing sequence of measurable cardinals and let $\mu=\bigcup_{\alpha\in\omega_1}\kappa_\alpha$.
Notice that $\mu>\cf(\mu)=\omega_1$.
We claim that $\binom{\mu}{\aleph_2}\rightarrow \binom{\mu}{\aleph_2}_2$.
To see this, suppose that $c:\mu\times\omega_2\rightarrow 2$ is a coloring.
Let $\mathscr{U}$ be an $\aleph_2$-complete ultrafilter over $\aleph_2$.
For every $\alpha\in\omega_1$ let $c_\alpha=c\upharpoonright(\kappa_\alpha\times\omega_2)$.
For every $\alpha\in\omega_1$ let $\mathscr{W}_\alpha$ be the ultrafilter generated by the $\omega$-closed unbounded subsets of $\kappa_\alpha$.
As in the first part of the proof, let $(A_\alpha,B_\alpha,i_\alpha)$ be so that $A_\alpha\in\mathscr{W}_\alpha, B_\alpha\in\mathscr{U}, i_\alpha\in\{0,1\}$ and $c_\alpha\upharpoonright(A_\alpha\times B_\alpha)$ is constantly $i_\alpha$.
Observe that the triple $(A_\alpha,B_\alpha,i_\alpha)$ is determined by the coloring $c$ and by the ultrafilters, so we do not use choice while creating the sequence $\langle(A_\alpha,B_\alpha,i_\alpha):\alpha\in\omega_1\rangle$.

Fix a set $I\subseteq\omega_1, |I|=\aleph_1$ and $i\in\{0,1\}$ such that $\alpha\in I\Rightarrow i_\alpha=i$.
Define $A=\bigcup\{A_\alpha:\alpha\in I\}$ and $B=\bigcap\{B_\alpha:\alpha\in I\}$.
Notice that $|A|=\mu$ and $B\in\mathscr{U}$ since $\mathscr{U}$ is $\aleph_2$-complete.
In particular, $|B|=\aleph_2$.
One can verify that $c''(A\times B)=\{i\}$, so the proof is accomplished.

\hfill \qedref{thmadcube}

The above theorem shows that $\binom{\mu}{\aleph_2}\rightarrow \binom{\mu}{\aleph_2}_2$ is independent over ZF for some $\mu>\cf(\mu)=\omega_1$, where the positive direction has been proved in models without choice and the negative direction holds in models of ZFC.
For example, the model of the generalized Martin's axiom with $2^{\omega_1}>\omega_2$ satisfies the negative relation.
We indicate, however, that the negative relation may hold even if the axiom of choice fails.
Indeed, all the relations of Theorem \ref{thmadcube} fail in Gitik's model \cite{MR576462}, as every uncountable cardinal has countable cofinality in this model.

The number of colors in the above theorem is optimal.
The trivial example of $c:\omega_2\times\omega_1\rightarrow\omega_1$ defined by $c(\alpha,\beta)=\beta$ shows that one cannot improve $2^{\aleph_0}$ to $\omega_1$ in the first statement.
Similarly, one cannot get infinitely many colors in the second part (though the statement holds with every finite number of colors).
Finally, by slicing $\mu>\cf(\mu)=\omega_1$ to $\omega_1$-many segments, each of which is of size less than $\mu$, one can show that $\binom{\mu}{\aleph_2}\nrightarrow \binom{\mu}{\aleph_2}_{\omega_1}$.

Ahead of dealing with the cube problem, we raise the possibility that the last part of the above theorem holds at $\mu=\aleph_{\omega_1}$.
It is not clear whether one can prove $\binom{\aleph_{\omega_1}}{\aleph_2}\rightarrow \binom{\aleph_{\omega_1}}{\aleph_2}_2$ under AD, but it seems that one can obtain this relation in ZF using the methods of \cite{MR1346109}.
For this end, one has to force $\aleph_{\alpha+1}$ to be measurable for every countable limit ordinal $\alpha\in\omega_1$, along with ${\rm DC}_{\aleph_1}$ and $\aleph_2$ being measurable.
Unlike AD in which we know who is the normal ultrafilter (so there is no need to choose it), if we force with the method of \cite{MR1346109} we must choose our measurable ultrafilter over each $\aleph_{\alpha+1}$, and hence ${\rm DC}_{\aleph_1}$ is needed.

\begin{question}
\label{qalephomega1} Is it consistent relative to ZF that $\binom{\aleph_{\omega_1}}{\aleph_2}\rightarrow \binom{\aleph_{\omega_1}}{\aleph_2}_2$?
\end{question}

We move now to the cube problem.
The cube version of the polarized relation is denoted by $\left( \begin{smallmatrix} \alpha \\ \beta \\ \gamma \end{smallmatrix} \right) \rightarrow \left( \begin{smallmatrix} \delta \\ \varepsilon \\ \zeta \end{smallmatrix} \right)_\theta$ and stipulates that for every coloring $c:\alpha\times\beta\times\gamma\rightarrow\theta$ there are $A\subseteq\alpha,B\subseteq\beta$ and $C\subseteq\gamma$ such that ${\rm otp}(A)=\delta,{\rm otp}(B)=\varepsilon,{\rm otp}(C)=\gamma$ and $c\upharpoonright(A\times B\times C)$ is constant.

Positive cube relations seem harder to achieve than the parallel relations applied to pairs.
In particular, it is unknown whether $\left( \begin{smallmatrix} \aleph_1 \\ \aleph_1 \\ \aleph_1 \end{smallmatrix} \right) \rightarrow \left( \begin{smallmatrix} \aleph_0 \\ \aleph_0 \\ \aleph_0 \end{smallmatrix} \right)_2$ is consistent.
In this case, it is also unknown whether the corresponding negative $\left( \begin{smallmatrix} \aleph_1 \\ \aleph_1 \\ \aleph_1 \end{smallmatrix} \right) \nrightarrow \left( \begin{smallmatrix} \aleph_0 \\ \aleph_0 \\ \aleph_0 \end{smallmatrix} \right)_2$ is consistent.
This problem is labeled as Question 28 in \cite{MR0280381}, and it appears in the generalized form as to whether $\left( \begin{smallmatrix} \kappa^+ \\ \kappa^+ \\ \kappa^+ \end{smallmatrix} \right) \rightarrow \left( \begin{smallmatrix} \kappa \\ \kappa \\ \kappa \end{smallmatrix} \right)_2$ in \cite[p. 110]{MR3075383}.
As noted by one of the referees, the status of this Question in ZFC is somewhat vague. It is claimed in several places that the answer is negative in ZFC, e.g. \cite{MR1297180}, \cite{MR648300} and \cite[Chapter 4]{MR3242261}.
However, no formal proof is available.

We shall try to prove positive instances of the cube relation under AD.
More precisely, we need the fact that $\kappa^+$ is a measurable cardinal and this holds at many places under determinacy assumptions, though it can also be forced without determinacy.

\begin{theorem}
\label{thmcube} if $\kappa^+$ is measurable then $\left( \begin{smallmatrix} \kappa^+ \\ \kappa^+ \\ \kappa^+ \end{smallmatrix} \right) \rightarrow \left( \begin{smallmatrix} \kappa \\ \kappa \\ \kappa \end{smallmatrix} \right)_2$. \newline
Moreover, $\left( \begin{smallmatrix} \kappa^+ \\ \kappa^+ \\ \kappa^+ \end{smallmatrix} \right) \rightarrow \left( \begin{smallmatrix} \kappa \\ \kappa \\ \kappa^+ \end{smallmatrix} \right)_2$.
\end{theorem}

\par\noindent\emph{Proof}. \newline
We prove the second statement, which of course implies the first one.
For this end, let $\mathscr{U}$ be a $\kappa^+$-complete ultrafilter over $\kappa^+$.
We define $\mathscr{U}^2$ over $\kappa^+\times\kappa^+$ as follows.
If $S\subseteq\kappa^+\times\kappa^+$ and $\beta\in\kappa^+$ then we let $S_\beta=\{\gamma\in\kappa^+:(\beta,\gamma)\in S\}$.
Now we define $S\in\mathscr{U}^2$ iff $\{\beta\in\kappa^+:S_\beta\in\mathscr{U}\}\in\mathscr{U}$.
It is routine to check that $\mathscr{U}^2$ is an ultrafilter and it is $\kappa^+$-complete.

Suppose that $c:\kappa^+\times\kappa^+\times\kappa^+\rightarrow 2$ is a coloring.
For every $\alpha\in\kappa^+,i\in\{0,1\}$ let $S^i_\alpha=\{(\beta,\gamma):c(\alpha,\beta,\gamma)=i\}$.
Notice that $(S^0_\alpha,S^1_\alpha)$ is a partition of $\kappa^+\times\kappa^+$ for every $\alpha\in\kappa^+$.
Hence for each $\alpha\in\kappa^+$ there is a unique $i(\alpha)\in\{0,1\}$ such that $S^{i(\alpha)}_\alpha\in\mathscr{U}^2$.
Consequently, there is a set $A\subseteq\kappa^+,|A|=\kappa^+$ and there is a fixed color $i\in\{0,1\}$ such that $\alpha\in A\Rightarrow i(\alpha)=i$.

Choose $a\in[A]^\kappa$.
For every $\alpha\in a$ let $B_{\alpha i}=\{\beta\in\kappa^+:(S^i_\alpha)_\beta\in\mathscr{U}\}$, so $B_{\alpha i}\in\mathscr{U}$.
Define $B=\bigcap\{B_{\alpha i}:\alpha\in a\}$ and conclude that $B\in\mathscr{U}$ since $\mathscr{U}$ is $\kappa^+$-complete.
Choose $b\in[B]^\kappa$ and define $E=\bigcap\{(S^i_\alpha)_\beta:\alpha\in a, \beta\in b\}$.
Again, $E\in\mathscr{U}$ so in particular $|E|=\kappa^+$.

Let us show that $c''(a\times b\times E)=\{i\}$.
Indeed, if $\alpha\in a,\beta\in b$ and $\gamma\in E$ then $\gamma\in(S^i_\alpha)_\beta$.
By the definition of these sets we see that $(\beta,\gamma)\in S^i_\alpha$.
By the definition of the coloring we see that $c(\alpha,\beta,\gamma)=i$, so the proof is accomplished.

\hfill \qedref{thmcube}

We can derive the following:

\begin{corollary}
\label{corad} Assume AD. \newline
Then $\left( \begin{smallmatrix} \aleph_1 \\ \aleph_1 \\ \aleph_1 \end{smallmatrix} \right) \rightarrow \left( \begin{smallmatrix} \aleph_0 \\ \aleph_0 \\ \aleph_1 \end{smallmatrix} \right)_2$ and $\left( \begin{smallmatrix} \aleph_{\omega+1} \\ \aleph_{\omega+1} \\ \aleph_{\omega+1} \end{smallmatrix} \right) \rightarrow \left( \begin{smallmatrix} \aleph_\omega \\ \aleph_\omega \\ \aleph_{\omega+1} \end{smallmatrix} \right)_2$.
\end{corollary}

\hfill \qedref{corad}

Remark that similar statements can be proved for $n$-copies of the same measurable cardinal where $3<n\in\omega$, using $\mathscr{U}^{n-1}$ instead of $\mathscr{U}^2$.
Likewise, the order type of the first two components can be a bigger ordinal since we can choose $a,b$ such that ${\rm otp}(a),{\rm otp}(b)\geq\delta$ for every $\delta\in\kappa^+$.
So actually we have, under AD and in these parameters, the optimal cube relation $\left( \begin{smallmatrix} \aleph_1 \\ \aleph_1 \\ \aleph_1 \end{smallmatrix} \right) \rightarrow \left( \begin{smallmatrix} \tau \\ \tau \\ \aleph_1 \end{smallmatrix} \right)_2$ for every $\tau\in\omega_1$.

Naturally, there are other interesting questions in the vicinity of the above results. Questions with infinite exponent were studies under AD in our context, see \cite{MR2984052}.
One of the referees posed an interesting question which lurks right around the corner, and we conclude the paper with an answer which shows the discrepancy between AD and weak fragments of choice.

\begin{theorem}
\label{thmrefereeq} Under AD we have $\left( \begin{smallmatrix} \omega_1 \\ \omega_1 \\ \omega \end{smallmatrix} \right) \rightarrow \left( \begin{smallmatrix} \omega_1 \\ \tau \\ \omega \end{smallmatrix} \right)_m$ for every $m\in\omega$ and every $\tau\in\omega_1$, but if ${\rm AC}_{\aleph_1}$ holds then $\left( \begin{smallmatrix} \omega_1 \\ \omega_1 \\ \omega \end{smallmatrix} \right) \nrightarrow \left( \begin{smallmatrix} \omega_1 \\ \omega \\ \omega \end{smallmatrix} \right)_2$.
\end{theorem}

\par\noindent\emph{Proof}. \newline
For the first statement all we need is just the measurability of $\aleph_1$, so under AD we can use Solovay's theorem.
Let $\mathscr{U}$ be an $\aleph_1$-complete ultrafilter over $\aleph_1$.
For simplicity we assume that $m=2$ and we indicate that the proof works at any finite number of colors.

Fix a countable ordinal $\tau\in\omega_1$.
For every $n\in\omega$ let $c_n=c\upharpoonright\omega_1\times\omega_1\times\{n\}$, so one can think of $c_n$ as a function from $\omega_1\times\omega_1$ into $\{0,1\}$.
For every $n\in\omega$ and for every $\beta\in\omega_1$ there is $i^n_\beta\in\{0,1\}$ such that $A^n_\beta=\{\alpha\in\omega_1:c_n(\alpha,\beta)=i^n_\beta\}\in\mathscr{U}$.
Hence for each $n\in\omega$ there is a set $B_n\in\mathscr{U}$ and some $i_n\in\{0,1\}$ for which $\beta\in B_n\Rightarrow i^n_\beta=i_n$.

For $m\neq n$ it is possible that $i_m\neq i_n$, but one can find a set $C\in[\omega]^\omega$ and a fixed $i\in\{0,1\}$ so that $n\in C\Rightarrow i_n=i$.
Let $B=\bigcap_{n\in\omega}B_n$, so $B\in\mathscr{U}$ by $\aleph_1$-completeness.
Let $B_\tau$ be the set which consists of the first $\tau$ elements of $B$.
Define $A=\bigcap\{A^n_\beta:n\in C,\beta\in B_\tau\}$ and notice again that $A\in\mathscr{U}$ since $\mathscr{U}$ is $\aleph_1$-complete.

If $\alpha\in A,\beta\in B_\tau$ and $n\in C$ then $c(\alpha,\beta,n)=c_n(\alpha,\beta)=i^n_\beta$ since $A\subseteq A^n_\beta$.
But $i^n_\beta=i_n$ since $\beta\in B_\tau\subseteq B\subseteq B_n$ and $i_n=i$ since $n\in C$.
Therefore, $c\upharpoonright(A\times B\times C)$ is constantly $i$, as desired.

To appreciate the above result suppose that ${\rm AC}_{\aleph_1}$ holds and choose a one-to-one function $g_\tau:\tau\rightarrow\omega$ for every $\tau\in\omega_1$.
We define $c:\omega_1\times\omega_1\times\omega\rightarrow 2$ as follows.
Given $\alpha,\beta\in\omega_1$ we let $\tau=\tau_{\alpha,\beta}=\max\{\alpha,\beta\}$ and $\sigma=\sigma_{\alpha,\beta}=\min\{\alpha,\beta\}$.
Now if $\alpha,\beta\in\omega_1$ and $n\in\omega$ we define $c(\alpha,\beta,n)=0$ iff $\sigma<\tau$ and $g_\tau(\sigma)\leq n$.

Suppose that $A\in[\omega_1]^{\omega_1},B\in[\omega_1]^\omega$ and $C\in[\omega]^\omega$.
Choose $\alpha\in A$ and $\beta\in B$ such that $\alpha>\beta$.
Let $\ell=g_\alpha(\beta)$.
Since $C$ is unbounded in $\omega$ one can find $n\in C$ so that $n\geq\ell$ and then $c(\alpha,\beta,n)=0$.
Now fix $\alpha\in A$ such that $\alpha>\sup(B)$.
For every $\beta\in B$ we have $\tau_{\alpha,\beta}=\alpha$ and hence $\sigma_{\alpha,\beta}=\beta$.
Fix $n\in C$.
Since $g_\alpha:\alpha\rightarrow\omega$ is one-to-one, the set $\{g_\alpha(\beta):\beta\in B\}$ is an infinite subset of $\omega$.
Hence one can find $\beta\in B$ for which $g_\alpha(\beta)>n$ and then $c(\alpha,\beta,n)=1$.
The negative relation is, therefore, established.

\hfill \qedref{thmrefereeq}

\newpage

\bibliographystyle{amsplain}
\bibliography{arlist}

\providecommand{\bysame}{\leavevmode\hbox to3em{\hrulefill}\thinspace}
\providecommand{\MR}{\relax\ifhmode\unskip\space\fi MR }
\providecommand{\MRhref}[2]{%
  \href{http://www.ams.org/mathscinet-getitem?mr=#1}{#2}
}
\providecommand{\href}[2]{#2}
\begin{thebibliography}{10}

\bibitem{MR2768693}
Uri Abraham and Menachem Magidor, \emph{Cardinal arithmetic}, Handbook of set
  theory. {V}ols. 1, 2, 3, Springer, Dordrecht, 2010, pp.~1149--1227.
  \MR{2768693}

\bibitem{MR2984052}
Arthur~W. Apter, Stephen~C. Jackson, and Benedikt L\"{o}we, \emph{Cofinality
  and measurability of the first three uncountable cardinals}, Trans. Amer.
  Math. Soc. \textbf{365} (2013), no.~1, 59--98. \MR{2984052}

\bibitem{MR1346109}
Arthur~W. Apter and Menachem Magidor, \emph{Instances of dependent choice and
  the measurability of {$\aleph_{\omega+1}$}}, Ann. Pure Appl. Logic
  \textbf{74} (1995), no.~3, 203--219. \MR{1346109}

\bibitem{MR1034564}
James~E. Baumgartner, \emph{Polarized partition relations and almost-disjoint
  functions}, Logic, methodology and philosophy of science, {VIII} ({M}oscow,
  1987), Stud. Logic Found. Math., vol. 126, North-Holland, Amsterdam, 1989,
  pp.~213--222. \MR{1034564}

\bibitem{MR1024901}
Hans-Dieter Donder and Jean-Pierre Levinski, \emph{Some principles related to
  {C}hang's conjecture}, Ann. Pure Appl. Logic \textbf{45} (1989), no.~1,
  39--101. \MR{1024901}

\bibitem{MR0280381}
P.~Erd\H{o}s and A.~Hajnal, \emph{Unsolved problems in set theory}, Axiomatic
  {S}et {T}heory ({P}roc. {S}ympos. {P}ure {M}ath., {V}ol. {XIII}, {P}art {I},
  {U}niv. {C}alifornia, {L}os {A}ngeles, {C}alif., 1967), Amer. Math. Soc.,
  Providence, R.I., 1971, pp.~17--48. \MR{0280381}

\bibitem{MR648300}
Paul Erd\H{o}s, \emph{Problems and results on finite and infinite combinatorial
  analysis. {II}}, Logic and algorithmic ({Z}urich, 1980), Monograph. Enseign.
  Math., vol.~30, Univ. Gen\`eve, Geneva, 1982, pp.~131--144. \MR{648300}

\bibitem{MR0202613}
P.~Erd{\H{o}}s, A.~Hajnal, and R.~Rado, \emph{Partition relations for cardinal
  numbers}, Acta Math. Acad. Sci. Hungar. \textbf{16} (1965), 93--196.
  \MR{MR0202613 (34 \#2475)}

\bibitem{MR795592}
Paul Erd{\H{o}}s, Andr{\'a}s Hajnal, Attila M{\'a}t{\'e}, and Richard Rado,
  \emph{Combinatorial set theory: partition relations for cardinals}, Studies
  in Logic and the Foundations of Mathematics, vol. 106, North-Holland
  Publishing Co., Amsterdam, 1984. \MR{MR795592 (87g:04002)}

\bibitem{MR3610266}
Shimon Garti, \emph{Amenable colorings}, Eur. J. Math. \textbf{3} (2017),
  no.~1, 77--86. \MR{3610266}

\bibitem{MR3201820}
Shimon Garti and Saharon Shelah, \emph{Partition calculus and cardinal
  invariants}, J. Math. Soc. Japan \textbf{66} (2014), no.~2, 425--434.
  \MR{3201820}

\bibitem{MR3509813}
\bysame, \emph{Open and solved problems concerning polarized partition
  relations}, Fund. Math. \textbf{234} (2016), no.~1, 1--14. \MR{3509813}

\bibitem{MR576462}
M.~Gitik, \emph{All uncountable cardinals can be singular}, Israel J. Math.
  \textbf{35} (1980), no.~1-2, 61--88. \MR{576462}

\bibitem{MR0371652}
R.~Laver, \emph{Partition relations for uncountable cardinals {$\leq 2^{\aleph
  _{0}}$}}, Infinite and finite sets ({C}olloq., {K}eszthely, 1973; dedicated
  to {P}. {E}rd{\H o}s on his 60th birthday), {V}ol. {II}, North-Holland,
  Amsterdam, 1975, pp.~1029--1042. Colloq. Math. Soc. Jan\'os Bolyai, Vol. 10.
  \MR{0371652 (51 \#7870)}

\bibitem{MR673792}
Richard Laver, \emph{An {$(\aleph _{2},\,\aleph _{2},\,\aleph _{0})$}-saturated
  ideal on {$\omega _{1}$}}, Logic {C}olloquium '80 ({P}rague, 1980), Stud.
  Logic Foundations Math., vol. 108, North-Holland, Amsterdam-New York, 1982,
  pp.~173--180. \MR{673792}

\bibitem{MR3242261}
R.~Daniel Mauldin (ed.), \emph{The {S}cottish {B}ook}, second ed.,
  Birkh\"{a}user/Springer, Cham, 2015, Mathematics from the Scottish Caf\'{e}
  with selected problems from the new Scottish Book, Including selected papers
  presented at the Scottish Book Conference held at North Texas University,
  Denton, TX, May 1979. \MR{3242261}

\bibitem{MR0297577}
Karel Prikry, \emph{On a problem of {E}rd{\H{o}}s, {H}ajnal and {R}ado},
  Discrete Math. \textbf{2} (1972), 51--59. \MR{0297577}

\bibitem{MR0505492}
S.~Shelah, \emph{A weak generalization of {MA} to higher cardinals}, Israel J.
  Math. \textbf{30} (1978), no.~4, 297--306. \MR{0505492 (58 \#21606)}

\bibitem{MR1318912}
Saharon Shelah, \emph{Cardinal arithmetic}, Oxford Logic Guides, vol.~29, The
  Clarendon Press, Oxford University Press, New York, 1994, Oxford Science
  Publications. \MR{1318912}

\bibitem{MR2768698}
John~R. Steel, \emph{An outline of inner model theory}, Handbook of set theory.
  {V}ols. 1, 2, 3, Springer, Dordrecht, 2010, pp.~1595--1684. \MR{2768698}

\bibitem{MR1278025}
Franklin~D. Tall, \emph{Some applications of a generalized {M}artin's axiom},
  Topology Appl. \textbf{57} (1994), no.~2-3, 215--248. \MR{1278025
  (95m:03101)}

\bibitem{MR1297180}
Stevo Todor\v{c}evi\'{c}, \emph{Some partitions of three-dimensional
  combinatorial cubes}, J. Combin. Theory Ser. A \textbf{68} (1994), no.~2,
  410--437. \MR{1297180}

\bibitem{MR3075383}
Neil~H. Williams, \emph{Combinatorial set theory}, Studies in Logic and the
  Foundations of Mathematics, vol.~91, North-Holland Publishing Co., Amsterdam,
  1977. \MR{3075383}

\end{thebibliography}

\end{document}